\documentclass[11pt,reqno]{amsart}
\usepackage{amsmath,amsopn,amssymb,amsthm,multicol  } 
\usepackage{hyperref}
\usepackage{color}
\usepackage{graphicx}

\newcommand{\fr}{\mathfrak}
\newcommand{\op}{\operatorname}

\DeclareMathOperator{\SO}{SO}
\DeclareMathOperator{\s}{S}
\DeclareMathOperator{\OO}{O}
\DeclareMathOperator{\Sp}{Sp}
\DeclareMathOperator{\SU}{SU}
\DeclareMathOperator{\U}{U}

 \newtheorem{lemma} {Lemma} [section]
\newtheorem{theorem}[lemma]{Theorem} 
\newtheorem{remark}[lemma] {Remark} 
\newtheorem{prop} [lemma]{Proposition}  
 
\newtheorem{corol}[lemma] {Corollary} 
 
\newtheorem{problem}[lemma] {Problem}

\begin{document}

\title[A review of compact geodesic orbit manifolds]{A review of compact geodesic orbit manifolds and the g.o. condition for $\SU(5)/\s(\U(2)\times \U(2))$}

\author{Andreas Arvanitoyeorgos, Nikolaos Panagiotis Souris and Marina Statha}
\address{University of Patras, Department of Mathematics, GR-26500 Rion, Greece, and 
Hellenic Open University, Aristotelous 18, GR-26335 Patras, Greece}
\email{arvanito@math.upatras.gr}
\address{University of Patras, Department of Mathematics, GR-26500 Rion, Greece}
\email{nsouris@upatras.gr  }
\address{University of Patras, Department of Mathematics, GR-26500 Rion, and
University of Thessaly,  Department of Mathematics, GR-35100 Lamia, Greece}
\email{statha@math.upatras.gr} 
\medskip

\begin{abstract}  Geodesic orbit manifolds (or g.o. manifolds) are those Riemannian manifolds $(M,g)$ whose geodesics are integral curves of Killing vector fields. Equivalently, there exists a Lie group $G$ of isometries of $(M,g)$ such that any geodesic $\gamma$ has the simple form $\gamma(t)=e^{tX}\cdot p$, where $e$ denotes the exponential map on $G$.  The classification of g.o. manifolds is a longstanding problem in Riemannian geometry. In this brief survey, we present some recent results and open questions on the subject focusing on the compact case. In addition we study the geodesic orbit condition for the space $\SU(5)/\s(\U(2)\times \U(2))$.


\end{abstract}

\maketitle
 

\section{Introduction}
\markboth{Andreas Arvanitoyeorgos Nikolaos Panagiotis Souris and Marina Statha}{A review of compact geodesic orbit manifolds and the g.o. condition for $SU(5)/S(U(2)\times U(2))$}

The prime examples of g.o. manifolds are the Euclidean space and the standard sphere. In fact, any connected g.o. manifold $(M,g)$ is \emph{homogeneous}, i.e. there exists a Lie group $G$ of isometries of $(M,g)$ acting transitively on $M$.  Then any geodesic of $(M,g)$ admits the form 

\begin{equation*}\gamma(t)=\exp(tX)\cdot p,\end{equation*}

\noindent where $\exp$ denotes the Lie exponential map on $G$, $p=\gamma(0)\in M$ and $\cdot$ denotes the (isometric) action of $G$ on $M$. The group $G$ can be the full isometry group of $(M,g)$ or possibly a proper Lie subgroup of the isometry group. In any case, $M$ is diffeomorphic to the coset $G/H$ where $H$ is the stabilizer of a point $p\in M$, i.e. the closed subgroup of $G$ whose elements fix the point $p\in M$. The corresponding Riemannian space $(M=G/H,g)$ is called a \emph{geodesic orbit space} while the ($G$-invariant) metric $g$ is called a \emph{geodesic orbit metric}.

Various well studied classes of homogeneous Riemannian manifolds are geodesic orbit, including \emph{symmetric spaces}, \emph{weakly symmetric spaces}, \emph{isotropy irreducible spaces}, \emph{$\delta$-homogeneous spaces} and \emph{naturally reductive spaces}. While most of the aforementioned classes have been completely determined, the complete classification of Riemannian g.o. manifolds remains a longstanding open problem.

The systematic study of g.o. manifolds was initiated by Kowalski and Vanhecke in 1991 (\cite{Ko-Va}) who classified the g.o. spaces up to dimension six.  Up to this day g.o. manifolds have been extensively studied (see for example the recent book \cite{Be-Ni-3} on the subject and the references therein). Since every g.o. manifold $(M,g)$ is essentially a homogeneous g.o. space $(G/H,g)$, the classification of g.o. manifolds reduces to the classification of g.o. spaces.

In this paper we focus on the compact case.  In Section \ref{conditions}, we discuss the classification problem of compact g.o. spaces in a more technical manner and we give some useful necessary and sufficient conditions for g.o. metrics.  In Section \ref{overview}, we give an overview of some of the progress so far on the classification of compact g.o. spaces.  In Section \ref{statement}, we state some of our recent results on the subject.  Finally in Section \ref{example} we study the geodesic orbit metric condition for  the homogeneous space $\SU(5)/\s(\U(2)\times \U(2))$.  

\medskip
\noindent 
{\bf Acknowledgment.}  This research is co-financed by Greece and the European Union (European Social Fund- ESF) through the Operational Programme ``Human Resources Development, Education and Lifelong Learning 2014-2020" in the context of the project ``Geodesic orbit metrics on homogeneous spaces of classical Lie groups" (MIS 5047124). 

\section{Compact g.o. spaces: Necessary and sufficient conditions}\label{conditions}  

We recall that a Riemannian homogeneous space $(G/H,g)$ is a Riemannian homogeneous manifold $M=G/H$ equipped with a \emph{$G$-invariant metric} $g$, i.e. a Riemannian metric that is invariant by the action of $G$ on $G/H$.  When studying compact g.o. spaces $(G/H,g)$, one has to firstly consider the following general problem.

\begin{problem}\label{problem1}Let $G/H$ be a homogeneous space with $G$ compact.  Find all $G$-invariant g.o. metrics $g$ on $G/H$.\end{problem}

Several partial classification results on g.o. spaces involve the solution of Problem \ref{problem1} under general assumptions for $G$ and/or $H$ (c.f. Section \ref{overview}).  We proceed to discuss some technical aspects of the problem.

Firstly, let $\fr{g}$ and $\fr{h}$ be the Lie algebras of $G$ and $H$ respectively. Since $G$ is compact, there exists an $\op{Ad}$-invariant inner product $Q$ on $\fr{g}$, which we consider fixed (here $\op{Ad}:G\rightarrow \op{Aut}(G)$) denotes the adjoint representation).  We consider a $Q$-orthogonal decomposition

\begin{equation*}\fr{g}=\fr{h}\oplus \fr{m}_{Q}=\fr{h}\oplus \fr{m},\end{equation*}

\noindent where $\fr{m}=\fr{m}_{Q}$ is the $Q$-orthogonal complement of $\fr{h}$ in $\fr{g}$ and coincides with the tangent space of $G/H$ at the origin $eH$.  Moreover, $\op{Ad}_H\fr{m}\subseteq \fr{m}$ and $[\fr{h},\fr{m}]\subseteq \fr{m}$.  The $G$-invariant metrics $g$ on $G/H$ are in bijection with $\op{Ad}_H$-invariant inner products $\langle \ ,\ \rangle$ on $\fr{m}$, which in turn are in bijection with \emph{metric endomorphisms} $\Lambda=\Lambda_Q\in \op{End}(\fr{m})$ satisfying 

\begin{equation}\label{metend}\langle X,Y\rangle=Q(\Lambda X, Y) \ \ \makebox{for all} \ \ X,Y \in \fr{m}.\end{equation} 
 
\noindent Each metric endomorphism $\Lambda$ is symmetric with respect to $Q$, positive definite and $\operatorname{Ad}_H$-equivariant (i.e. $\Lambda$ commutes with $\operatorname{Ad}_h$ for all $h\in H$).  Conversely, any endomorphism of $\fr{m}$ with the above properties corresponds to a unique $G$-invariant Riemannian metric on $G/H$ through Equation \eqref{metend}. For the rest of this paper, we will make no distinction between a $G$-invariant metric on $G/H$ and its corresponding metric endomorphism $\Lambda \in\op{End}(\fr{m})$.  The following is a necessary and sufficient condition for a metric endomorphism to define a g.o. metric.

\begin{prop}\label{GOCond}\emph{(\cite{Al-Ar}, \cite{So1})} Let $G$ be a compact Lie group.  The Riemannian space $(G/H,\Lambda)$ is a g.o. space if and only if there exists a map $\xi:\fr{m}\rightarrow \fr{h}$ such that 

\begin{equation}\label{GO}[\xi(X)+X,\Lambda X]=0 \ \ \makebox{for all} \ \ X\in \fr{m}.\end{equation}
\end{prop}  

\noindent Therefore, Problem \ref{problem1} reduces to the following.

\begin{problem}\label{problem2}Let $G/H$ be a homogeneous space with $G$ compact.  Find all metric endomorphisms $\Lambda\in \op{End}(\fr{m})$ satisfying equation \eqref{GO}.\end{problem}

From equation \eqref{GO}, we observe that there are two fundamental difficulties in solving Problem \ref{problem2}:

\noindent \textbf{1.} The general form of the metric endomorphism $\Lambda$ may be quite complicated in higher dimensions (see for example \cite{So1}), which makes the solution of Equation \eqref{GO} impossible through straightforward means.\\
\noindent \textbf{2.} The map $\xi:\fr{m}\rightarrow \fr{h}$, which is called a \emph{geodesic graph}, depends on the embedding of $H$ in $G$ and it is generally non-differentiable at the origin. As a result, even if one could solve Equation \eqref{GO} through straightforward or computational means, it is not possible to do so for large classes of spaces simultaneously.

To remedy difficulty \textbf{1.}, various simplification techniques from geometry, Lie theory and representation theory have been applied, e.g. \emph{polar representations} (\cite{Ta}), \emph{principal orbit types} (\cite{CheNiNi}), \emph{root theory of semisimple Lie algebras} (\cite{Al-Ar}), the \emph{isotypic decomposition} (\cite{So1}) and several others (see for example \cite{Be-Ni-3}, \cite{Ni2}, \cite{So1} and references therein).  These methods aim to create strong necessary conditions for a metric endomorphism $\Lambda$ to be a g.o. metric. Although we will not delve into much detail about those methods, we will state one of the most important simplification conditions called the \emph{normalizer lemma}.

\begin{lemma}\label{NormalizerLemma}\emph{(\cite{Ni2})} 
The inner product $\langle \ ,\ \rangle$  generating the metric of a geodesic orbit Riemannian space $(G/H,g)$, is not only $\op{Ad}(H)$-invariant but also $\op{Ad}(N_G(H^0))$-invariant, where $N_G(H^0)$ is the normalizer of the identity component $H^0$ of the group $H$ in $G$.\end{lemma}

As a result of the above lemma, we have the following (\cite{Ar-Sou-St2}). 

\begin{prop}\label{conclusion} {\rm Decompose the tangent space $\fr{m}$ of $G/H$ into the $Q$-orthogonal sum $\fr{m}=\fr{n}\oplus \fr{p}$, where $\fr{p}$ is the tangent space of $G/N_G(H^0)$ and $\fr{n}$ is the Lie algebra of the group  $N_G(H^0)/H^0$.  Then any metric endomorphism $\Lambda$, corresponding to a g.o. metric on $G/H$, has the block-diagonal form
$$\Lambda=\begin{pmatrix}\left.\Lambda\right|_{\fr{n}} & 0\\
0& \left.\Lambda\right|_{\fr{p}}\end{pmatrix},
$$}

\noindent where the block $\left.\Lambda\right|_{\fr{n}}$ defines a bi-invariant metric on the group $N_G(H^0)/H^0$ and the block $\left.\Lambda\right|_{\fr{p}}$ defines a $G$-invariant g.o. metric on the space $G/N_G(H^0)$.
\end{prop}
 
An application of Proposition \ref{conclusion} will be shown in Section \ref{example}.\\

The following lemma is also useful, since it allows us to prove that some eigenvalues of a g.o. metric are equal under certain algebraic conditions.

\begin{lemma}\label{EigenEq}\emph{(\cite{So1})}
 Let $(G/H,g)$ be a g.o. space with $G$ compact and with corresponding metric endomorphism $\Lambda$ with respect to an $\op{Ad}$-invariant inner product $Q$. Let $\fr{m}$ be the $Q$-orthogonal complement of $\fr{h}$ in $\fr{g}$.\\ 
\textbf{1.} Assume that $\fr{m}_1,\fr{m}_2$ are $\operatorname{ad}(\fr{h})$-invariant, pairwise $Q$-orthogonal subspaces of $\fr{m}$ such that $[\fr{m}_1,\fr{m}_2]$ has non-zero projection on ${(\fr{m}_1\oplus \fr{m}_2)^\bot}$. Let $\lambda_1,\lambda_2$ be eigenvalues of $\Lambda$ such that $\left.\Lambda\right|_{\fr{m}_i}=\lambda_i\op{Id}$, $i=1,2$.  Then $\lambda_1=\lambda_2$.\\
\textbf{2.} Assume that $\fr{m}_1,\fr{m}_2,\fr{m}_3$ are $\operatorname{ad}(\fr{h})$-invariant, pairwise $Q$-orthogonal subspaces of $\fr{m}$ such that $[\fr{m}_1,\fr{m}_2]$ has non-zero projection on $\fr{m}_3$. Let $\lambda_1,\lambda_2,\lambda_3$ be eigenvalues of $\Lambda$ such that $\left.\Lambda\right|_{\fr{m}_i}=\lambda_i\op{Id}$, $i=1,2,3$.  Then $\lambda_1=\lambda_2=\lambda_3$.
\end{lemma}

To address difficulty \textbf{2.}, the usual strategy involves studying Problem \ref{problem2} for suitably chosen families of spaces $G/H$, such that the isotropy groups $H$ have similar embeddings in $G$ or have consistent properties (e.g. $H$ is abelian).  This approach has led to various partial classifications (see Section \ref{overview}), which have significantly advanced our understanding of compact g.o. spaces and furthered the progress of their classification.

For any compact homogeneous space $G/H$ there exists at least one solution to Problem \ref{problem2}:  The metric $\Lambda_Q=\op{Id}$ satisfies Equation \eqref{GO}, with the corresponding geodesic graph $\xi$ being the zero map. A $G$-invariant metric with $\Lambda=\op{Id}$ is called a \emph{normal metric} and it is clearly a g.o. metric.  More generally, a $G$-invariant metric is called \emph{naturally reductive} if it satisfies Equation \eqref{GO} with $\xi$ linear.  The corresponding space $(G/H,\Lambda)$ is called a naturally reductive space.  It is evident that any naturally reductive metric is a g.o. metric.  Naturally reductive spaces have been extensively studied by geometers, while their classification is also an open problem (\cite{Ag}, \cite{St}).

\section{Overview of results on compact g.o. spaces}\label{overview}

Since the initial investigation of O. Kowalski and L. Vanhecke in \cite{Ko-Va},  several authors have contributed in the subject with quite interesting results. We will try to give a short summary for some of these below.
In \cite{Go} C. Gordon reduced the classification of g.o. Riemannian manifolds $G/H$  to the following cases:
(a) $G$ nilpotent, (b) $G$ non compact semisimple and (c) $G$ compact semisimple 
In \cite{Ta} H. Tamaru
 classified g.o. spaces fibered over irreducible symmetric spaces.
In \cite{DuKoNi} Z. Du\v sek, O. Kowalski and S. Nik\v cevi\'c gave examples of g.o. manifolds in dimension 7.
In \cite{Be-Ni-1} V. Berestovskii and  Yu. Nikonorov  investigated the g.o. property for genelarized normal homogeneous Riemannian manifolds (in other terminology $\delta$-homogeneous manifolds).
Also, in \cite{Al-Ar} D. Alekseevsky and  the first author 
classified  simply connected generalized flag manifolds admitting non normal g.o. metrics, obtaining two infinite families.

 Homogeneous geodesics in Heisenberg groups and other pseudo-Riemannian manifolds were studied by 
 Z. Du\v sek and O. Kowalski in \cite{DK}.
In \cite{CheNi}  Z. Chen and Yu. Nikonorov classified compact, simply connected g.o. spaces with two isotropy summands.
In \cite{GoNi} Gordon and Nikonorov gave a
 Geometric and algebraic characterization of g.o. manifolds that are diffeomorphic to $\mathbb{R}^n$. 
In \cite{CheCheDe}  H. Chen, Z. Chen and S. Deng and in \cite{Ni3} Nikonorov gave examples of left-invariant Einstein metrics on compact simple Lie groups which are not g.o.
Note that there are examples of homogeneous Einstein metrics that are neither naturally reductive, nor g.o.  (e.g. $\SU(3)/T_{\rm max}$, or Aloff-Wallach spaces $\SU(3)/\mathbb{S}^2_{k,l}$).

In \cite{ArvWa} the first author and Y. Wang
 classified g.o. spaces among generalized Wallach spaces.  Recall that these are homogeneous spaces $G/H$ such that $\mathfrak{m}=\mathfrak{m}_1\oplus\mathfrak{m}_2\oplus\mathfrak{m}_3$ with $[\mathfrak{m}_i, \mathfrak{m}_i]\subset\mathfrak{k}$, $\mathfrak{m}_i$ irreducible.
In \cite{ArvWaZa} the first author, Y. Wang and G. Zhao classified g.o. spaces among $M$-spaces.  These are homogeneous spaces $G/K_1$ such that $G/K$ is a generalized flag manifold with $K=C(S)\times K_1$, ($S$ a torus in $G$ and $K_1$ the semisimple part of $K$).
In \cite{So1} the second author gave classification of g.o. spaces
 by using  isotypic decomposition.
 
In \cite{NikNik} Yu. Nikolayevsky and Yu. Nikonorov proved that a Ledger-Obata space is a g.o. space if and only if it is naturally reductive.  This is a homogeneous space of the form $(F\times F\times\cdots\times F)/{\rm diag}(F)$, where $F$ is a connected, compact, simple Lie group.
In \cite{So2} the second author studied g.o. spaces $(G/H,g)$ with $H$ abelian, and in \cite{CheNiNi}
Z. Chen, Yu. Nikolayevski and Yu. Nikonorov
studied g.o. spaces $(G/H,g)$ with $H$ simple.

In \cite{CheCheZ} H. Chen, Z. Chen and  F. Zhu constructed g.o. spaces from by strongly isotropy spaces.

Finally, the notion of homogeneous geodesics has been extended to geodesics which are orbits of a product of two or more exponential factors, i.e. $\gamma (t)={\rm exp} (tX){\rm exp}tY\cdot o$ by  the first two authors and G. Calvaruso (\cite{ArCaSou}).

For $G$ compact and semisimple we have the following
recent classification results:

In \cite{So2} the second author proved the following: 
Let $G$ be a compact, semisimple Lie group and $H$ a closed \textbf{abelian} subgroup.  The Riemannian space $(G/H,\Lambda)$ is geodesic orbit if and only if $\Lambda$ is a normal metric, i.e. induced from a bi-invariant metric on $G$.

Also, in \cite{CheNiNi} Z. Chen, Yu. Nikolayevski and Yu Nikonorov proved the following: Let $(G/H,g)$ be a compact, irreducible and non normal g.o. space with $H$ \textbf{simple}.  Then $G/H$ is (up to a finite cover) one of the following spaces:

\begin{eqnarray*}\SO(9)/{\rm Spin}(7)& \SO(10)/{\rm Spin}(7) & \SO(11)/{\rm Spin}(7)\\
{\rm E}_6/{\rm Spin}(10) & \SU(n+p)/\SU(p)  & \SO(2n+1)/\SU(n)\\
\SO(4n + 2)/\SU(2n + 1) & \Sp(n + 1)/\Sp(n) & \SU(2n + 1)/\Sp(n) \\
 {\rm Spin}(8)/{\rm G}_2 & \SO(9)/{\rm G}_2. \end{eqnarray*}
 
 Therefore, the above lead to the following

\medskip
\underline{Open Question:}
Classify the g.o. spaces $(G/H,g)$ with $H$ \textbf{semisimple}.

\section{Statement of main results}\label{statement}

 The  following results have appeared in \cite{Ar-Sou-St} and \cite{Ar-Sou-St2}. 

\begin{theorem}\label{main1}
Let $G/H$ be the space $\SO(n)/\SO(n_1)\times \cdots \times \SO(n_s)$, where $0<n_1+\cdots +n_s\leq n$, and $n_j>1, j=1,\dots , s$.  A $G$-invariant Riemannian metric on $G/H$ is geodesic orbit if and only if it is a normal metric, i.e. it is induced from an $\op{Ad}$-invariant inner product on the Lie algebra $\fr{so}(n)$ of $\SO(n)$.
\end{theorem}

\begin{remark} {\rm We remark that if $n\neq 4$ then $\fr{so}(n)$ is simple, and thus any $\op{Ad}$-invariant inner product is homothetic to the negative of the Killing form $B(X,Y)=(n-2)\op{Trace}(XY)$. If $n=4$ then $\fr{so}(n)\equiv\fr{so}(3)\oplus \fr{so}(3)$, and thus any $\op{Ad}$-invariant inner product is homothetic to the negative of the one-parameter family $B_1+\lambda B_2$, $\lambda>0$, where $B_1$ denotes the Killing form of the first simple factor $\fr{so}(3)$ and $B_2$ denotes the Killing form of the second simple factor $\fr{so}(3)$.}
\end{remark}

As a result of Theorem \ref{main1} 
we obtain the following.

\begin{corol}\label{corol_o(n)}
Let $G/H$ be one of the spaces $\OO(n)/\OO(n_1)\times \cdots \times \OO(n_s)$ or $\SO(n)/\op{S}(\OO(n_1)\times \cdots \times \OO(n_s))$, where $0<n_1+\cdots +n_s\leq n$, $n_j>1$.  A $G$-invariant Riemannian metric on $G/H$ is geodesic orbit if and only if it is normal.
\end{corol} 

The second main result is the following.

\begin{theorem}\label{main2}Let $G/H$ be the space $\U(n)/\U(n_1)\times \cdots \times \U(n_s)$, where $n_1+\cdots+n_s\leq n$, and let $N_G(H)$ be the normalizer of $H$ in $G$. If $n_1+\dots +n_s=n$, then a $G$-invariant Riemannian metric on $G/H$ is geodesic orbit if and only if it is the normal metric induced from the $\op{Ad}$-invariant inner product $B(X,Y)=-\op{Trace}(XY)$ in $\fr{u}(n)$.  If $n_1+\cdots +n_s<n$, then a $G$-invariant Riemannian metric $g$ on $G/H$ is geodesic orbit if and only if $g=g_{\mu}$, $\mu>0$, where $g_{\mu}$ denotes a one-parameter family of deformations of the normal metric induced from the inner product $B$, along the center of the group $N_G(H)/H$.
\end{theorem}

\smallskip
We note that the metrics in Theorem \ref{main2} generalize the g.o. metrics on the Berger spheres $\U(n)/\U(n-1)$ (\cite{Ni1}) and the g.o. metrics on the complex Stiefel manifolds $\U(n)/\U(n-k)$ (\cite{So1}).  We also note that the g.o. metrics on the related class of real flag manifolds were recently studied in \cite{Neg}.
 Among other results, it is shown in \cite{Neg} that every g.o. metric on the real flag manifold $\SO(n)/\op{S}(\OO(n_1)\times \cdots \times \OO(n_s))$, $n_1+\cdots +n_s=n$, is normal, which is a special case of Corollary \ref{corol_o(n)}.

\begin{theorem}\label{main3}
Let $G/H$  be the space $\Sp(n)/\Sp(n_1)\times \cdots \times \Sp(n_s)$, where $0<n_1+\cdots +n_s\leq n$.  If $G/H\neq \Sp(n)/\Sp(n-1)$ (i.e. if $n-(n_1+\cdots+n_s)\neq 1$ or $s>1$) then a $G$-invariant Riemannian metric on $G/H$ is geodesic orbit if and only if it is the standard metric induced from the Killing form on the Lie algebra $\fr{sp}(n)$ of $\Sp(n)$.  

If $G/H=\Sp(n)/\Sp(n-1)$ (i.e. $s=1$ and $n-n_1=1$) then a $G$-invariant metric $g$ on $G/H$ is geodesic orbit if and only if $g=g_{\mu}$, $\mu>0$, where $g_{\mu}$ denotes a one-parameter family of deformations of the standard metric $g_1$, along the fiber $\Sp(1)$ of the fibration $\Sp(n)/\Sp(n-1)\rightarrow \Sp(n)/\Sp(1)\times \Sp(n-1)$.
\end{theorem}

We remark that the non standard geodesic orbit metric $g_{\mu}$ appears in  \cite{Nik0} and \cite{Ta}.

\section{Geodesic orbit metrics on the homogeneous space $\SU(5)/\s(\U(2)\times \U(2))$}\label{example}

We consider the homogeneous space $G/H=\SU(5)/\s(\U(2)\times \U(2))$.  Here the isotropy subgroup $H=\s(\U(2)\times \U(2))=\{A\in \U(2)\times \U(2):\det(A)=1\}$ is diagonally embedded in $G=\SU(5)$. We will study the $G=SU(5)$-invariant g.o. metrics on $G/H$. 

Let $\fr{g}=\fr{su}(5)$ be the Lie algebra of $G$ and let $\fr{h}=\{X\in \fr{u}(2)\times \fr{u}(2):\op{trace}(X)=0\}$ be the Lie algebra of $H$.  Moreover, we consider the negative $Q$ of the Killing form of $\fr{g}$, given by 

\begin{equation*}Q(X,Y)=-\op{trace}(XY).\end{equation*}

\noindent We also consider the $Q$-orthogonal decomposition 

\begin{equation*} \fr{g}=\fr{h}\oplus \fr{m}.\end{equation*}

 We need to explicitly describe the tangent space $\fr{m}$ of $G/H$ at the origin.  To this end, we firstly consider a $Q$-orthogonal basis of $\fr{g}$. Let $M_{5}\mathbb C$ be the set of complex $5\times 5$ matrices and let $E_{ab}\in M_{5}\mathbb C$, $a,b=1,\dots,5$, be the matrix having  1 in the $(a,b)$-entry and zero elsewhere.  For $a,b=1,\dots,5$, we set
\begin{equation}\label{mel2}e_{ab}=E_{ab}-E_{ba},\quad f_{ab}=\sqrt{-1}(E_{ab}+E_{ba}), \ \ a<b, \ \ f_{aa}=\frac{1}{2}\sqrt{-1}E_{aa}.
\end{equation}

Observe that
$e_{ab}=-e_{ba}, {f}_{ab}={f}_{ba}$.
 The set 
\begin{equation}\label{set2}\mathcal{B}=\left\{ e_{ab}, {f}_{cd}, f_{ll}-f_{l+1,l+1}:1\le a<b\le 5 ,\  1\le c < d\le 5, \ 1\le l \le 4 \right\},
\end{equation}  
 is a basis of $\fr{su}(5)$ which is orthogonal with respect to $Q$.  The non zero bracket relations among the vectors \eqref{mel2} are given by
 
\begin{equation}\label{cn}
\begin{array}{ccc}
 [e_{ab},e_{cd}]=\delta_{bc}e_{ad}-\delta_{ad}e_{cb}-\delta_{ac}e_{bd}-\delta_{bd}e_{ac},&
[{f}_{ab},e_{cd}]=\delta_{bc}{f}_{ad}-\delta_{ad}{f}_{cb}+\delta_{ac}{f}_{bc}-\delta_{bd}f_{ac},\\

\ \ \ [{f}_{ab}, {f}_{cd}]=-\delta_{bc}e_{ad}+\delta_{ad}e_{cb}-\delta_{ac}e_{bd}-\delta_{bd}e_{ac}.\end{array}
\end{equation}

We may identify the Lie algebra $\fr{h}$ with the $Q$-orthogonal sum

\begin{equation*}\fr{h}=\op{span}_{\mathbb R}\{f_{22}-f_{33},e_{23},f_{23}, f_{44}-f_{55}, e_{45},f_{45}\}\oplus \op{span}_{\mathbb R}\{v\}, \ \ v:=f_{22}+f_{33}-f_{44}-f_{55}.\end{equation*}

Using relations \eqref{cn}, one can verify that the normalizer $\fr{n}_{\fr{g}}(\fr{h})$ of $\fr{h}$ in $\fr{g}$ is the Lie algebra  $\fr{h}\oplus \op{span}_{\mathbb R}\{w\}$, where 

\begin{equation*}w:=4f_{11}-f_{22}-f_{33}-f_{44}-f_{55}.\end{equation*}

Therefore, we have the decomposition 

\begin{equation*}\fr{m}=\fr{n}\oplus \fr{p}, \ \ \makebox{where} \ \ \fr{n}=\op{span}_{\mathbb R}\{w\}\end{equation*}

\noindent coincides with the Lie algebra of the group $N_G(H^0)/H^0$ (c.f. Proposition \ref{conclusion}).  

Let $\Lambda$ be a g.o. metric on $G/H$.  From Proposition \ref{conclusion} and the fact that $\fr{n}$ is one dimensional, we have 

\begin{equation}\label{ln}\left.\Lambda\right|_{\fr{n}}=\mu\op{Id}.\end{equation}

 To fully describe the tangent space $\fr{m}$ with respect to the basis $\mathcal{B}$, it remains to find a basis for the space $\fr{p}$.  The space $\fr{p}$ decomposes into $\op{ad}_{\fr{n}_{\fr{g}}(\fr{h})}$-irreducible submodules as

\begin{equation*}\fr{p}=\fr{m}_{01}\oplus \fr{m}_{02}\oplus \fr{m}_{12},\end{equation*}

where the subspaces $\fr{m}_{ij}$ have the following description in terms of the basis $\mathcal{B}$:

\begin{eqnarray*}\fr{m}_{01}&=&\op{span}_{\mathbb R}\{e_{12},f_{12},e_{13},f_{13}\}, \\
\fr{m}_{02}&=&\op{span}_{\mathbb R}\{e_{14},f_{14},e_{15},f_{15}\},\\
\fr{m}_{12}&=&\op{span}_{\mathbb R}\{e_{24},f_{24},e_{25},f_{25},e_{34},f_{34},e_{35},f_{35}\}.\end{eqnarray*}

All the submodules $\fr{m}_{ij}$ are pairwise inequivalent with respect to $\op{ad}_{\fr{n}_{\fr{g}}(\fr{h})}$, and thus $\left.\Lambda\right|_{\fr{m}_{ij}}=\lambda_{ij}\op{Id}$.  By combining Proposition \ref{conclusion} with Lemma \ref{EigenEq}, we deduce that

\begin{equation}\label{lp}\left.\Lambda\right|_{\fr{p}}=\lambda\op{Id}.\end{equation}

Equations \eqref{ln} and \eqref{lp} imply that up to homothety we have

\begin{equation}\label{L}\Lambda=\begin{pmatrix}\mu\left.\op{Id}\right|_{\fr{n}} & 0\\
0 & \left.\op{Id}\right|_{\fr{p}}\end{pmatrix}.\end{equation}

Now for a vector of the form $X=\sum_{1\leq i<j\leq 5}(c_{ij}e_{ij}+d_{ij}f_{ij})\in \fr{g}$, denote by $\bar{X}$ the vector 

\begin{equation*}\bar{X}:=\sum_{1\leq i<j\leq 5}(-d_{ij}e_{ij}+c_{ij}f_{ij}).\end{equation*}

By virtue of equations \ref{cn}, we obtain the following.

\begin{lemma}\label{brackets}Let $X_{01}\in \fr{m}_{01}$, $X_{02}\in \fr{m}_{02}$, $X_{12}\in \fr{m}_{12}$. The following Lie-bracket relations are true.

\begin{eqnarray*}&&[w,\fr{h}]=\{0\}, \ \ \ \ \ \ \ [v,\fr{h}]=\{0\}, \ \ \ \ \ \ \ \ [v,w]=0\\
&& [w,X_{01}]=5\bar{X_{01}}  \ \ \ [w,X_{02}]=5\bar{X_{02}}, \ \ \ [w,X_{12}]=0,\\
&& [v,X_{01}]=-\bar{X_{01}}  \ \ \ [v,X_{02}]=\bar{X_{02}}, \ \ \ \ [v,X_{12}]=2\bar{X_{12}}.\end{eqnarray*}\end{lemma}

We are ready to find equivalent conditions for the g.o. metric $\Lambda$.  Let $X\in \fr{m}$.  According to the decomposition 

\begin{equation*}\fr{m}=\fr{n}\oplus \fr{p}=\op{span}_{\mathbb R}\{w\}\oplus \fr{m}_{01}\oplus \fr{m}_{02}\oplus\fr{m}_{12},\end{equation*}

\noindent write $X=cw+X_{01}+X_{02}+X_{12}$.  Without any loss of generality, we may assume that $c=1$ and thus $X=w+X_{01}+X_{02}+X_{12}$.  By Equation \eqref{L}, we have $\Lambda X=\mu w+X_{01}+X_{02}+X_{12}$.  By Proposition \ref{GOCond}, $\Lambda$ is a g.o. metric if and only if for all $X\in \fr{m}$ there exists $a\in \fr{h}$ such that $0=[a+X,\Lambda X]$.  By using the facts that $X=w+X_{01}+X_{02}+X_{12}$, $\Lambda X=\mu w+X_{01}+X_{02}+X_{12}$ and Lemma \ref{brackets}, Equation \eqref{GO} is equivalent to

\begin{eqnarray*}0&=&[a+X,\Lambda X]=[a+ w+X_{01}+X_{02}+X_{12},\mu w + X_{01}+X_{02}+X_{12}]\\
&=&[a,X_{01}]+5(1-\mu)\bar{X_{01}}+[a,X_{02}]+5(1-\mu)\bar{X_{02}}+[a,X_{12}].\end{eqnarray*}

Since $[a,X_{01}]\in \fr{m}_{01}$, $[a,X_{02}]\in \fr{m}_{02}$ and $[a,X_{12}]\in \fr{m}_{12}$, the above equation is equivalent to the system

\begin{eqnarray}1. [a,X_{01}]+5(1-\mu)\bar{X_{01}}&=&0,\\
   2.  [a,X_{02}]+5(1-\mu)\bar{X_{02}}&=&0,\\
3. [a,X_{12}]&=&0,
\end{eqnarray}

\noindent with unknown $a\in \fr{h}$.

\end{document}